\documentclass[11pt]{article}
\usepackage{mathrsfs}
\usepackage{latexsym}
\usepackage{amsmath,amsthm}
\usepackage{amsfonts}
\usepackage{url}
\usepackage{amssymb}
\usepackage{mytheo}
\def\longto{\mathop{\longrightarrow}\limits}
\def\clap#1{\hbox to 0pt{\hss#1\hss}}

\def\mathclap{\mathpalette\mathclapinternal}

\def\mathclapinternal#1#2{%
\clap{$\mathsurround=0pt#1{#2}$}}
\newcommand{\narsum}[1]{\sum_{\mathclap{#1}}}

\newcommand{\etab}{\bar\eta}

\usepackage{bbm}
\newcommand{\chr}{\boldsymbol{\mathbbm{1}}} %
\newcommand{\pred}[1]{\chr_{\left\{ #1 \right\}}}

\newcommand{\E}{\mathbb{E}}
\renewcommand{\P}{\mathbb{P}}
\newcommand{\var}{\mathrm{Var}}

\newcommand{\inv}{^{-1}} %

\newcommand{\TV}[1]{\nrm{#1}_{\textrm{{\tiny \textup{TV}}}}}

\newcommand{\Lip}[1]{\nrm{#1}_{\textrm{{\tiny \textup{Lip}}}}}

\usepackage{geometry}
\makeatletter
\newcommand{\ben}{\begin{enumerate}}
\newcommand{\een}{\end{enumerate}}
\newcommand{\bit}{\begin{itemize}}
\newcommand{\eit}{\end{itemize}}

\newcommand{\X}{\N}
\newcommand{\supr}[1]{^{(#1)}}
\newcommand{\seq}[3]{(#1_{#2},\ldots,#1_{#3})}
\newcommand{\sseq}[3]{#1_{#2}^{#3}}  %

\newcommand{\dsabs}[1]{\bigl| #1 \bigr|}
\newcommand{\nrm}[1]{\left\Vert #1 \right\Vert}

\newcommand{\calL}{\mathcal{L}}

\newcommand{\tha}{\theta}
\newcommand{\R}{\mathbb{R}}
\newcommand{\N}{\mathbb{N}}

\newcommand{\beq}{\begin{eqnarray*}}
\newcommand{\eeq}{\end{eqnarray*}}
\newcommand{\beqn}{\begin{eqnarray}}
\newcommand{\eeqn}{\end{eqnarray}}
\newcommand{\paren}[1]{\left( #1 \right)}
\newcommand{\sqprn}[1]{\left[ #1 \right]}
\newcommand{\tlprn}[1]{\left\{ #1 \right\}}
\newcommand{\set}[1]{\tlprn{#1}}
\newcommand{\abs}[1]{\left| #1 \right|}

\newcommand{\gn}{\, | \,}

\newcommand{\ts}{\textstyle}

\renewcommand{\th}{\ensuremath{^{\mathrm{th}}}~}

\newcommand{\hide}[1]{}
\newcommand{\oo}[1]{\frac{1}{#1}}
\def\eps{\varepsilon}
\newcommand{\bp}{\boldsymbol{p}}

\newcommand{\ninf}{\longto_{n\to\infty}}
\newcommand{\bx}{\boldsymbol{x}}
\newcommand{\by}{\boldsymbol{y}}
\title{
Uniform Chernoff and
Dvoretzky-Kiefer-Wolfowitz-type inequalities
for Markov chains and related processes}
\author{
Aryeh Kontorovich\footnote{
Department of Computer Science, Ben-Gurion University, Beer Sheva, Israel
}
\thanks{
This research was supported by the ISRAEL SCIENCE FOUNDATION (grant No. 1141/12).
} \and Roi Weiss$^*$}
\begin{document}
\maketitle
\begin{abstract}
We observe that the technique of Markov contraction
can 
be used to establish measure concentration for a broad class of
non-contracting chains. 
In particular, geometric ergodicity provides a simple and versatile framework.
This 
leads to a short, elementary
proof of a general concentration inequality for Markov and 
hidden Markov chains (HMM),
which supercedes some of the known results and
easily extends to other processes such as Markov trees.
As applications, we give a 
Dvoretzky-Kiefer-Wolfowitz-type inequality
and
a uniform Chernoff bound.
All of our bounds are dimension-free and hold for countably infinite state spaces.
\end{abstract}
\section{Introduction}
\subsection{Background}
The last decade or so has seen a flurry of activity in concentration of measure for non-independent processes.
A recent survey may be found in 
\cite{kontorovich12}, 
with pointers to more specialized surveys therein. 
Rather than recapitulating these surveys here, we shall proceed directly to 
the relevant
recent developments.
Let $X_1,X_2,\ldots$ be a sequence of $\N$-valued random variables obeying some joint law (distribution). 
Using the shorthand 
$\calL(\sseq{X}{j}{n}\gn \sseq{X}{1}{i}=x)$ to denote the law 
of $\seq{X}{j}{n}$ conditioned on $\seq{X}{1}{i}=x
\in\X^i$,
let us define, for
$n\in\N$,
$1\leq i<j\leq n$,
$y\in\X^{i-1}$ and 
$w,w'\in\X$, 
\beq
\eta_{ij}(y,w,w') &=&
\TV{
\calL(\sseq{X}{j}{n}\gn \sseq{X}{1}{i}={y w})-
\calL(\sseq{X}{j}{n}\gn \sseq{X}{1}{i}={y w'})
},
\eeq
(where $\TV{\cdot}=\oo2\nrm{\cdot}_1$ is the total variation norm) and
\beqn
\label{eq:etadef}
\etab_{ij} &=&
\sup_{y\in\X^{i-1},w,w'\in\X}
\eta_{ij}(y,w,w').
\eeqn
The coefficients $\etab_{ij}$, 
termed {\em $\eta$-mixing coefficients}
in \cite{kontram06},
play a central role in several recent concentration results. 
Define 
$\Delta$ to be 
the upper-triangular $n\times n$ matrix,
with $
\Delta_{ii}=1$
 and
$\Delta_{ij} = \etab_{ij}$
for $1\leq i<j\leq n$.
In 2007, 
\cite{chazottes07} and \cite{kontram06}
independently proved that
for any 
$f:\X^n\to\R$ with $\Lip{f}\leq 1$ with respect to the 
Hamming metric\footnote{
Meaning: if $x,y\in\N^n$ differ in only 1 coordinate then $\abs{f(x)-f(y)}\le1$,
see Section~\ref{sec:aux}.}, we have
\beqn
\label{eq:kontram}
P\paren{\abs{f-\E f}>n\eps } &\leq& 2\exp\paren{-\frac{2n\eps^2}{
\min\set{\nrm{\Delta}_2,\nrm{\Delta}_\infty}^2
}},
\eeqn
where $\nrm{\Delta}_p$ is the $\ell_p$ operator norm
(\cite{chazottes07} achieve
the better constant in the exponent, given here). Earlier, Samson 
\cite{samson00} had given a concentration result for convex $\ell_2$-Lipschiz
functions $f:[0,1]^n\to\R$, which likewise involved the coefficients $\etab_{ij}$,
and these are also implicit in Marton's earlier work
\cite{marton98,marton03,marton04}.
In order to apply (\ref{eq:kontram}) in a Markov setting, 
one must upper-bound 
$\nrm{\Delta}_2$
or
$\nrm{\Delta}_\infty$
for the Markov chain in question. 
The earliest such results relied on contraction. Let
$p(\cdot\gn\cdot)$ 
be the transition kernel associated with a given Markov chain,
and define the 
(D\"oblin)
{\em contraction coefficient}
\beqn
\label{eq:kapdef}
\kappa = 
\sup_{x,x'\in\X}\TV{p(\cdot\gn x)-p(\cdot\gn x')}.
\eeqn
It is shown in \cite{kontram06} and 
\cite{samson00} that $\etab_{ij}\le\kappa^{j-i}$ and therefore 
$\nrm{\Delta}_\infty\le(1-\kappa)\inv$; this implies the concentration
bound 
\beq
\P(\abs{f-\E f}>n\eps)\le2\exp(-2(1-\kappa)^2n\eps^2)
\eeq
for $1$-Lipschitz functions $f$,
which Marton \cite{marton96}
had essentially obtained earlier by other means.
The contraction method was pushed further to obtain concentration
results for hidden Markov chains \cite{kontram06},
undirected Markov chains and 
Markov tree processes \cite{kontorovich12}, but its
applicability requires the rather stringent condition that
$\kappa<1$. Already in \cite{marton98}, Marton observed
that a significantly weaker mixing condition suffices,
and yields
tighter and more informative bounds.
Indeed, consider a Markov chain with stationary distribution $\pi$
and 
conditional
$s$\th step distribution $\calL(X_s\gn X_1=x)$,
and define the ``inverse mixing time''\footnote{
This terminology is non-standard.
}
\beqn
\label{eq:taudef}
\tau_s = \sup_{x\in\X}\TV{ \calL(X_s\gn X_1=x) - \pi}.
\eeqn
A simple calculation (Lemma \ref{lem:eta-tau}) shows that
$\etab_{ij}\le2\tau_{j-i}$, and thus
\beq
\nrm{\Delta}_\infty 
-1
=
\max_{1<i<n}\sum_{j=i+1}^n\etab_{ij}
\le
2\max_{1<i<n}\sum_{j=i+1}^n
\tau_{j-i}
.
\eeq
A rich body of work deals with bounding $\tau_s$ via spectral 
\cite{0726.60069},
Poincar\'e \cite{springerlink:10.1007/BF02214660},
log-Sobolev \cite{0867.60043}
and Lyapunov \cite{springerlink:10.1007/s00440-011-0373-4}
methods, among others (see the references in the works cited). From our perspective,
the {\em geometric ergodicity} condition allows for the simplest exposition while sacrificing the least
generality. A Markov chain is said to be geometrically ergodic with constants 
$1\le G<\infty$ and $0\le\tha<1$ if
\beqn
\label{eq:tauGtheta}
\tau_s \le G\tha^{s-1}, \qquad s=1,2,\ldots
.
\eeqn
Any finite ergodic Markov chain is geometrically ergodic,
and the dependence of $G,\tha$ on various structural properties of the chain in question
is the subject of 
a diverse and prolific literature
(including the references above). 
We also stress that the geometric ergodicity assumption is largely dictated by
expositional convenience, since any non-trivial bound on the inverse mixing time $\tau_s$
will yield straightforward analogues of our results.
In this paper, we explore some consequences of geometric ergodicity as pertaining to
concentration and statistical inference for Markov and hidden Markov chains.
We leverage two basic insights: (i) 
even though hidden Markov chains are a 
considerably richer class of processes
than Markov chains (there exist HMMs not 
realizable by any finite-order Markov chain),
for the purposes of measure concentration, 
the underlying Markov chain is all that matters
and (ii) geometric ergodicity, while significantly more general than contractivity,
yields essentially the same concentration bounds.
Another advantage of our approach is its elementary nature:
taking the bound in (\ref{eq:kontram}) as a given, nothing beyond basic linear algebra is used.
Given the recent interest in prediction and parameter inference for 
HMMs \cite{DBLP:journals/corr/abs-1203-0683,DBLP:conf/colt/HsuKZ09,MR2244426,Siddiqi10,DBLP:conf/icml/KontorovichNW13,NIPS2013_5007},
our result have potential to be applicable beyond the abstract setting studied here.
Furthermore, since concentration results for Markov chains extend easily for other Markov-type processes
(such as trees \cite{kontorovich12}), our results here should extend to those as well.
\subsection{Main results}
\paragraph{Concentration.}
Our first result is a concentration inequality for hidden Markov chains, which generalizes many of the previous such bounds.
We will henceforth write 
``$(G,\theta)$-geometrically ergodic'' as shorthand for
``geometrically ergodic with constants $1\le G<\infty$ and $0\le\tha<1$''.
Hidden Markov chains and their associated notions of stationarity and geometric ergodicity are formally 
defined in Section \ref{sec:notation}.
\bethn
\label{thm:main-conc}
Let $Y_1,Y_2,\ldots$ be a $\N$-valued hidden Markov chain whose underlying 
$\N$-valued
Markov chain is 
$(G,\theta)$-geometrically ergodic.
Then, for any $n\in\N$ and $f:\N^n\to\R$ with 
$\Lip{f}\le1$
(under the Hamming metric), 
we have
\beq
\label{eq:main-conc}
\P\paren{{f(\sseq{Y}{1}{n})-\E f(\sseq Y1n)}>n\eps} 
&\leq& \exp\paren{-\frac{
n(1-\tha)^2\eps^2
}{
2G^2
}},
\eeq
with an identical bound for the other tail.
\enthn
Although the result in Theorem \ref{thm:main-conc} does not appear to have been published anywhere,
it is a simple consequence of widely known facts (we give a proof in Section \ref{sec:methods} for completeness).
Our main contribution 
lies in the apparently novel applications.
\paragraph{DKW-type inequality.}
Let us recall the 
Dvoretzky-Kiefer-Wolfowitz inequality
\cite{MR0083864,MR1062069},
stated here for the discrete case.
Suppose $X_1,X_2,\ldots$ are iid $\N$-valued random variables with common distribution function $F$,
and define the empirical distribution function $\hat F_n$ induced by 
$\seq X1n$:
\beq
\hat F_n(x) = \oo n\sum_{i=1}^n \pred{X_i\le x},\qquad x\in\N.
\eeq
The DKW inequality states that
\beq
\P\paren{\sup_{x\in\N}\abs{\hat F_n(x)-F(x)}>\eps}\le2\exp(-2n\eps^2),\qquad\eps>0,n\in\N.
\eeq
We present the following 
Markovian version
of this inequality. 
\bethn
\label{thm:dkw}
Let $Y_1,Y_2,\ldots$ be 
a stationary 
$\N$-valued
$(G,\theta)$-geometrically ergodic 
Markov or hidden Markov 
chain 
with 
stationary distribution 
$\rho\in\R^\N$. 
For $n\in\N$, define $\hat\rho\supr{n}\in\R^\N$ to be the empirical estimate of $\rho$:
\beqn
\label{eq:rhodef}
\hat\rho\supr{n}_y=\oo n\sum_{i=1}^n \pred{Y_i=y},
\qquad y\in\N.
\eeqn
Then
\beq
\P\paren{
\nrm{\rho-\hat\rho\supr{n}}_\infty
>
\sqrt{\frac{1+2G\tha}{n(1-\tha)}}
+
\eps
}
\le \exp\paren{-\frac{
n(1-\tha)^2\eps^2
}{
2G^2
}},
\qquad
n\in\N,
\eps>0
.
\eeq
\enthn
Note that a naive application of Theorem \ref{thm:main-conc} to each 
$\hat\rho\supr{n}_y$ individually, 
combined with the union bound,
yields
\beqn
\label{eq:rho-sup-naive}
\P\paren{\nrm{\rho-\hat\rho\supr{n}}_\infty>\eps}\le 
2
\nrm{\rho}_0
\exp\paren{-\frac{
n(1-\tha)^2\eps^2
}{
2G^2}},
\eeqn
where $\nrm{\rho}_0$ is the 
number of non-zero entries in
$\rho$. 
The bound in (\ref{eq:rho-sup-naive}) is vacuous
for $\rho$ with infinite support.
The assumption that the chain starts
in the
stationary distribution is not at all restrictive, as shown in Section \ref{sec:stat}.
\paragraph{Uniform Chernoff bound.}
Let $Y_1,Y_2,\ldots$ be a stationary 
$\N$-valued 
$(G,\theta)$-geometrically ergodic 
Markov or hidden Markov chain as above, and 
consider the occupation frequency:
\beq
\hat \rho\supr{n}(E) = \oo n\sum_{i=1}^n \pred{Y_i\in E},
\qquad E\subseteq \N.
\eeq
A naive application of Theorem \ref{thm:main-conc} might yield a deviation bound along the lines of
\beq
\P\paren{\abs{
\rho(E)-\hat \rho\supr{n}(E)
}>\eps}\le 
2
|E|
\exp\paren{-\frac{
n(1-\tha)^2\eps^2
}{
2|E|^2G^2}},
\eeq
where $|E|$ is the cardinality of $E$ and $\rho$ is the stationary distribution as above.
We will give a much stronger bound, that is not only independent of $E$ but is actually uniform over all $E\subseteq\N$.
\bethn
\label{thm:u-chernoff}
Define
\beq
\Lambda_n(\rho) = 
\gamma_n(G,\tha)
\narsum{\rho_y\ge1/n}\sqrt{\rho_y}
+
\min\set{
\gamma_n(G,\tha)
\narsum{\rho_y<1/n}\sqrt{\rho_y}~,
\narsum{\rho_y<1/n}\rho_y}
,
\qquad n\in\N,
\eeq
where
\beq
\gamma_n(G,\tha) = \oo2\sqrt{\frac{1+2G\tha}{n(1-\tha)}}.
\eeq
Then:
\bit
\item[(a)]
for all distributions $\rho\in\R^\N$,
\beq
\lim_{n\to\infty} \Lambda_n(\rho)=0,
\eeq
\item[(b)]
\beq
\P\paren{
\sup_{E\subseteq\N}
\abs{
\rho(E)-\hat \rho\supr{n}(E)
}>
\Lambda_n(\rho)+
\eps}\le 
\exp\paren{-\frac{
n(1-\tha)^2\eps^2
}{
2G^2}}.
\eeq
\hide{}
\eit
\enthn
We remark that the rate at which $\Lambda_n(\rho)$ decays to $0$ depends on $\rho$ and
may be arbitrarily slow for heavy-tailed distributions. 
When
$\sum_{y\in\N}\sqrt{\rho_y}<\infty$,
we get a simpler estimate in (b) via
\beq
\Lambda_n(\rho) \le
\gamma_n(G,\tha)
\narsum{y\in\N}
\sqrt{\rho_y}
.
\eeq
Again, the stationarity assumption is quite mild (Section \ref{sec:stat}).
\subsection{Related work}
In parallel to the work on concentration of measure results for Markov chains 
\cite{MR2424985,adamczak12,MR2511280,kontram06,marton96,samson00},
grew a body of independent results on Chernoff-type bounds for these processes.
The papers \cite{0829.60022,Dinwoodie1998585,DBLP:journals/siamcomp/Gillman98,DBLP:journals/cpc/Kahale97,0938.60027}
played a founding role,
and various extensions and refinements followed \cite{1056.60070,DBLP:journals/cpc/Wagner08}.
In a remarkable recent development \cite{DBLP:conf/stacs/ChungLLM12}, 
optimal Chernoff-Hoeffding bounds are obtained based on the mixing time at a constant threshold.
Concentration of Lipschitz functions of mixing sequences, with applications to the 
Kolmogorov-Smirnov statistic, were considered in \cite{MR1771956}. The paper \cite{MR2759184}
examines the concentration of empirical distributions for non-independent sequences
satisfying 
Poincar\'e or
log-Sobolev inequalities.
\hide{}
\section{Methods and proofs}
\label{sec:methods}
\subsection{Preliminaries}
\label{sec:notation}
For readability, we will sometimes write the matrix entry $A_{x,y}$
as $A(x\gn y)$. We will use the terms {\em hidden Markov chain} and HMM interchangeably.
\paragraph{Markov chains.}
We will represent Markov kernels by column-stochastic $\N\times\N$ matrices denoted by the letter $A$.
Thus, a Markov chain with transition kernel $A$ and initial distribution $p_1$ induces the following
distribution on $\N^n$:
\beqn
\label{eq:markX}
\calL(X_1,\ldots,X_n) = p_1(X_1)\prod_{i=1}^{n-1} A({X_{i+1}\gn X_{i}}).
\eeqn
\paragraph{Hidden Markov chain.}
A hidden Markov chain (also known as hidden Markov model [HMM]) is specified by
the triple $(p_1,A,B)$, where $(p_1,A)$ are the Markov chain parameters as above
and $B$ is an $\N\times\N$ column-stochastic matrix of {\em emission probabilities}.
This HMM induces a distribution on $\N^n$ as follows. 
Let $X\in\N^n$ be distributed according to (\ref{eq:markX}) and define the conditional distribution
$\calL(\cdot\gn X)$ over $Y\in\N^n$:
\beq
\calL(Y\gn X) 
= 
\prod_{i=1}^{n} B({Y_i\gn X_i}).
\eeq
It follows that 
\beq
\calL(Y)=\sum_{x\in\N^n} \P(X=x)\calL(Y\gn X=x).
\eeq
We will refer to $Y$ as a {\em hidden Markov} chain and to $X$ as its {\em underlying} Markov chain.
\paragraph{Stationary distributions and chains.}
The stationary distribution $\pi\in\R^\N$ of the Markov chain 
with transition kernel $A$ is the unique stochastic vector
satisfying $A\pi=\pi$.
The Markov chain
induced by $(p_1,A)$ is said to be stationary if $p_1=\pi$. 
It is well-known that, for ergodic Markov chains,
\beq
\pi = \lim_{n\to\infty}\calL(X_n) = \lim_{n\to\infty} \E\hat\pi\supr{n},
\eeq
where
\beq
\hat\pi\supr{n}_x=\oo n\sum_{i=1}^n \pred{X_i=x},
\qquad x\in\N.
\eeq
In the geometrically ergodic case,
observing that
$\E\hat\pi\supr{n}=\oo n\sum_{i=1}^n \calL(X_i)$, we have
\beq
\TV{\E\hat\pi\supr{n}-\pi} &=& 
\TV{\oo n\sum_{i=1}^{n}(\calL(X_i)-\pi)}\\
&\le& 
\oo n\sum_{i=1}^{n}\TV{\calL(X_i)-\pi} \\
&=&
\oo n\sum_{i=1}^{n}\TV{\sum_{x\in\N}\calL(X_i\gn X_1=x)p_1(x)-\pi} \\
&\le&
\oo n\sum_{i=1}^{n}\sum_{x\in\N}p_1(x)\TV{\calL(X_i\gn X_1=x)-\pi} \\
&\le&
\oo n\sum_{i=1}^{n}\sum_{x\in\N}p_1(x)G\tha^{i-1} 
=
\frac{G}{(1-\theta)n}.
\eeq
For a hidden Markov chain, we define the stationary distribution
$\rho=B\pi$, and observe that
\beq
\rho = \lim_{n\to\infty}\calL(Y_n) = \lim_{n\to\infty} \E\hat\rho\supr{n},
\eeq
where $\hat\rho\supr{n}$ is defined in (\ref{eq:rhodef}).
Since $\hat\rho\supr{n}$ is distributed as $B\hat\pi\supr{n}$,
we have
\beqn
\label{eq:hmm-nonstat}
\TV{\E\hat\rho\supr{n}-\rho}
\le
\TV{\E\hat\pi\supr{n}-\pi}
\le \frac{G}{(1-\theta)n}.
\eeqn
The bound in
(\ref{eq:hmm-nonstat}) 
suggests that, at least to some degree,
the statistical behavior of an HMM is controlled by its underlying Markov chain.
We expand upon this observation further:
\belen
\label{lem:hmm<markov}
Let $X$ and $X'$ be two Markov chains induced by $(\xi,A)$ and $(\xi',A')$, respectively.
For a given emission matrix $B$, let $Y$ and $Y'$ be the hidden Markov chains induced by
$(\xi,A,B)$ and $(\xi',A',B)$.
Then
\beq
\TV{\calL(Y_{i\in I}) - \calL(Y'_{i\in I})}
\le
\TV{\calL(X_{i\in I}) - \calL(X'_{i\in I})} 
,
\qquad 
I\subseteq\set{1,\ldots,n},
n\in\N.
\eeq
\enlen
\bepf
Immediate from Jensen's inequality, since hidden Markov chains are
convex mixtures of Markov chains.
\enpf
The proofs of Theorems~\ref{thm:dkw} and \ref{thm:u-chernoff} will
require bounds on 
$\nrm{\hat\rho\supr{n}-\rho}$, but unlike in (\ref{eq:hmm-nonstat}),
the expectation is on the outside of the norm.
\subsection{Markov contraction}
\label{sec:contr}
Let us recast the 
contraction coefficient defined in
(\ref{eq:kapdef}) in the language of Markov kernels:
\beq
\kappa = 
\sup_{x,x'\in\N}\TV{A(\cdot\gn x)-A(\cdot\gn x')}.
\eeq
The term ``contraction'' is justified by the following 
simple fact \cite{bremaud99,kontram06}:
\belen
[Markov, 1906 \cite{mar1906}]
\label{lem:markov-contraction}
For any two stochastic vectors $\xi,\psi\in\R^\N$, we have
\beq
\TV{A(\xi-\psi)}\le\kappa\TV{\xi-\psi}.
\eeq
\enlen
Our principal application
of this result will be in the context of geometrically ergodic Markov kernels.
\becon
\label{cor:ergo-contr}
Let $A$ be a $(G,\tha)$-geometrically ergodic Markov kernel. Then for all $n\in\N$,
the $n$-step kernel $A^n$ has contraction coefficient $\kappa\le2G\tha^{n}$.
\encon
\bepf
Let $\pi$ be the stationary distribution of $A$
and $\xi,\psi\in\R^\N$ 
two point masses.
Then
\beq
\TV{A^n\xi-A^n\psi} &\le & \TV{A^n\xi-\pi}+\TV{A^n\psi-\pi} \\
&\le& 2\tau_{n+1}\le 2G\tha^{n}.
\eeq
\enpf
\subsection{Proof of main inequality}
\label{sec:main-conc}
In this section, we prove Theorem \ref{thm:main-conc}.
The first order of business is to bound the $\eta$-mixing coefficient by
the inverse mixing time,
and hence in terms of $G$ and $\tha$.
\belen
\label{lem:eta-tau}
Let $Y$ be a $(G,\tha)$-geometrically ergodic hidden Markov chain
and let $\etab_{ij}$ and $\tau_s$ be as defined in (\ref{eq:etadef}) and (\ref{eq:taudef}), respectively.
Then
\beq
\etab_{ij}\le2\tau_{j-i+1}\le 2G\tha^{j-i}
,
\qquad
n\in\N, 1\le i<j\le n.
\eeq
\enlen
\bepf
Let $X$ be the Markov chain underlying $Y$ and endow $\etab_{ij}(X)$, $\etab_{ij}(Y)$ 
with the obvious meaning.
Then \cite[Theorem 7.1]{kontram06} shows that
\beq
\etab_{ij}(Y)\le
\etab_{ij}(X).
\eeq
Next, 
Remark 4 and the Theorem preceding it in \cite{kontorovich12} show that
\beq
\etab_{ij}(X)\le \kappa(A^{j-i})
\eeq
where $\kappa(A^{j-i})$ is the contraction coefficient of the $(j-i)$-step Markov kernel of $X$.
Finally, 
Corollary \ref{cor:ergo-contr} yields
\beq
\kappa(A^{j-i}) \le 2\tau_{j-i+1} \le 2G\tha^{j-i}.
\eeq
\enpf
\bepf[Proof of Theorem \ref{thm:main-conc}.]
By (\ref{eq:kontram}), it suffices to upper-bound
\beq
\nrm{\Delta}_\infty 
=
1+\max_{1<i<n}\sum_{j=i+1}^n\etab_{ij}
.
\eeq
Applying Lemma \ref{lem:eta-tau}, we get
\beq
\max_{1<i<n}\sum_{j=i+1}^n\etab_{ij} &\le& 2G\max_{1<i<n}\sum_{j=i+1}^n \tha^{j-i} \\
&\le& 
2G\sum_{k=1}^\infty\tha^k.
\eeq
Since $G\ge1$ by assumption, we have
\beq
1+2G\sum_{k=1}^\infty\tha^k &\le& 2G\sum_{k=0}^\infty\tha^k \\
&\le& 
\frac{2G}{1-\tha}.
\eeq
\enpf
\subsection{Proof of the DKW-type inequality}
\newcommand{\xiy}{\xi\supr{y}}
In this section, we prove Theorem \ref{thm:dkw}.
Let $Y_1,Y_2,\ldots$ be a stationary
$(G,\theta)$-geometrically ergodic 
hidden Markov 
chain 
with 
stationary distribution $\rho$,
and define the $\set{0,1}$-indicator variables
\beqn
\label{eq:xidef}
\xiy_i = \pred{Y_i=y},
\qquad i,y\in\N.
\eeqn
Then $\hat\rho$, defined in (\ref{eq:rhodef}), is given by
$\hat\rho_y
=\oo n\sum_{i=1}^n\xiy_i$, where we have dropped the superscript $(n)$ from $\hat\rho$ for readability.
Observing that the map $\seq Y1n\mapsto n\nrm{\rho-\hat\rho}_\infty$ is $1$-Lipschitz under the Hamming metric (Lemma~\ref{lem:Flip}),
we apply Theorem \ref{thm:main-conc}:
\beq
\P(\nrm{\rho-\hat\rho}_\infty
>\E\nrm{\rho-\hat\rho}_\infty+\eps)
\le
\exp\paren{-\frac{n(1-\tha)^2\eps^2}{2G^2}}.
\eeq
Hence, it remains to bound $\E\nrm{\rho-\hat\rho}_\infty$.
\belen
\label{lem:suprho}
\beq
\E\nrm{\rho-\hat\rho}_\infty
&\le&
\sqrt{\frac{1+2G\tha}{n(1-\tha)}}.
\eeq
\enlen
{\em Remark.} This estimate is nearly optimal:
in the case where $Y_i$ are iid (i.e., $\tha=0$)
Bernoulli variables with parameter $p$, we have \cite[Theorem 1]{Berend2013}
\beq
\sqrt{\frac{p(1-p)}{2n}}
\le
\E\nrm{\rho-\hat\rho}_\infty 
\le
\sqrt{\frac{p(1-p)}{n}}
,
\qquad n\ge 2,~p\in[1/n,1-1/n].
\eeq
\bepf
Jensen's inequality yields
\beqn
\paren{\E\nrm{\rho-\hat\rho}_\infty}^2 
&\le&
\E\sqprn{ \nrm{\rho-\hat\rho}_\infty^2} \nonumber\\
&\le&
\E\sqprn{ \sum_{y\in\N}\abs{\rho_y-\hat\rho_y}^2} 
\label{eq:rho2}
\\
&=&
\nonumber
\sum_{y\in\N}\E\paren{\rho_y-\hat\rho_y}^2 
=
\sum_{y\in\N}\var[\hat\rho_y].
\eeqn
Putting $S\supr{y}_n=\sum_{i=1}^n\xiy_i$, we have
\beqn
\label{eq:varid}
n^2\var[\hat\rho_y] = \E{\paren{S\supr{y}_n}^2} - 
\paren{\E{S\supr{y}_n}}^2
\eeqn
and
\beqn
\label{eq:Eid}
\E{S\supr{y}_n}=n\rho_y.
\eeqn
To bound $
\E{\paren{S\supr{y}_n}^2}
$, we compute
\beqn
\E{\paren{S\supr{y}_n}^2}
&=& \E\sqprn{\sum_{1\le i,j\le n}\xiy_i\xiy_j} \nonumber\\
&=& \sum_{i=1}^n\E{\paren{\xiy_i}^2}+2\narsum{1\le i<j\le n}\E\sqprn{\xiy_i\xiy_j} \nonumber\\
&=& 
\label{eq:S2id}
n\rho_y
+2\narsum{1\le i<j\le n}\E\sqprn{\xiy_i\xiy_j},
\eeqn
where the last identity holds since $\xiy_i\in\set{0,1}$.
It now remains to estimate $\E\sqprn{\xiy_i\xiy_j}$.
To this end, we claim that
\beq
\nrm{\calL(Y_i\gn Y_1=y)-\rho}_\infty \le G\tha^{i-1},
\qquad
i,y\in\N.
\eeq
Indeed,
denoting the parameters of $Y$ by $(\pi,A,B)$ and letting
$X$ be the underlying Markov chain, we have
\beq
\nrm{\calL(Y_i\gn Y_1=y_1)-\rho}_\infty &\le&
\TV{\calL(Y_i\gn Y_1=y_1)-\rho}\\
&=&
\oo2\sum_{y_i\in\N}\abs{\P(Y_i=y_i\gn Y_1=y_1)-\rho_{y_i}} \\
&=&
\oo2\sum_{y_i\in\N}
\abs{
\sum_{x_i\in\N} 
B_{y_i,x_i}
\paren{
\P(X_i=x_i\gn Y_1=y_1)-\pi_{x_i}}
} \\
&\le&
\oo2\sum_{y_i\in\N}
\sum_{x_i\in\N} 
B_{y_i,x_i}
\abs{\P(X_i=x_i\gn Y_1=y_1)-\pi_{x_i}}\\
&=&
\oo2\sum_{x_i\in\N }\abs{\P(X_i=x_i\gn Y_1=y_1)-\pi_{x_i}}\\
&=&
\TV{\sum_{x_1\in\N}\calL(X_i\gn X_1=x_1)\P(X_1=x_1\gn Y_1=y_1)-\pi}\\
&\le&
\sup_{x_1\in\N}\TV{\calL(X_i\gn X_1=x_1)-\pi}
\le G\tha^{i-1}.
\eeq
Hence,
\beq
\E\sqprn{\xiy_i\xiy_j} &=& \P(Y_i=y,Y_j=y)\\
&=& \P(Y_1=y,Y_{j-i+1}=y)\\
&=& \P(Y_1=y)\P(Y_{j-i+1}=y\gn Y_1=y)\\
&\le& \rho_y(\rho_y+G\tha^{j-i}),
\eeq
and therefore
\beqn
\narsum{1\le i<j\le n}\E\sqprn{\xiy_i\xiy_j}
&=&
\sum_{k=1}^{n-1}(n-k)\P(Y_1=y)\P(Y_{k+1}=y\gn Y_1=y)\nonumber\\
&\le&
\sum_{k=1}^{n-1}(n-k)\rho_y(\rho_y+G\tha^{k}) \nonumber\\
&=& 
\frac{n(n-1)}{2}\rho_y^2
+
\frac{G\tha}{1-\tha}
\paren{n-\frac{1-\tha^n}{1-\tha}}\rho_y \nonumber\\
&\le&
\label{eq:xixi}
\frac{n(n-1)}{2}\rho_y^2
+n\frac{G\tha}{1-\tha}\rho_y.
\eeqn
Combining (\ref{eq:varid}), (\ref{eq:Eid}), (\ref{eq:S2id}), and (\ref{eq:xixi}), we have
\beq
\var[\hat\rho_y]&\le&
\oo{n^2}\paren{n\rho_y+n(n-1)\rho_y^2+2n\frac{G\tha}{1-\tha}\rho_y-n^2\rho_y^2}\\
&=& \frac{\rho_y}{n}\paren{1-\rho_y+\frac{2G\tha}{1-\tha}}\\
&\le& {\rho_y}\frac{1+2G\tha}{n(1-\tha)}.
\eeq
Since $\sum_{y\in\N}\rho_y=1$, the claim follows from (\ref{eq:rho2}).
\enpf
{\em Remark.} Note that in the process of proving a deviation estimate on $\nrm{\rho-\hat\rho}_\infty$,
we have actually proven a stronger one --- namely, for the $\ell_2$ norm.
\subsection{Proof of the uniform Chernoff bound}
In this section, we prove Theorem \ref{thm:u-chernoff}.
As before, $Y_1,Y_2,\ldots$ is a stationary
$(G,\theta)$-geometrically ergodic 
hidden Markov 
chain 
with 
stationary distribution $\rho$.
Since by Lemma~\ref{lem:Flip}
the map
$\seq Y1n\mapsto n\TV{\rho-\hat\rho}$ is $1$-Lipschitz under the Hamming metric,
Theorem \ref{thm:main-conc} applies:
\beqn
\label{eq:conc-tv}
\P(\TV{\rho-\hat\rho}
>\E\TV{\rho-\hat\rho}+\eps)
\le
\exp\paren{-\frac{n(1-\tha)^2\eps^2}{2G^2}}.
\eeqn
As before, the crux of the matter is to bound $\E\TV{\rho-\hat\rho}$.
Recall the definition of $\Lambda_n$ from the statement of Theorem \ref{thm:u-chernoff}.
\belen
\label{lem:rho-Lam}
\beq
\E\TV{\rho-\hat\rho}\le\Lambda_n.
\eeq
\enlen
{\em Remark.} This bound is nearly optimal:
when the $Y_i$ are iid,
we have \cite[Proposition 3]{Berend2013}
\beq
\E\TV{\rho-\hat\rho}
\ge
\oo4\Lambda_n-\oo{8\sqrt n}
,
\qquad n\ge 2,~p\in[1/n,1-1/n].
\eeq
\bepf
We proceed by breaking up the expectation into two terms,
\beqn
\label{eq:2terms}
\E\TV{\rho-\hat\rho} &=&
{\ts\oo2}\narsum{y:\rho_y<1/n} \E\abs{\rho_y-\hat\rho_y}
+
{\ts\oo2}\narsum{y:\rho_y\ge1/n} \E\abs{\rho_y-\hat\rho_y},
\eeqn
and bounding each term separately.
To bound the second term, we note,
as in the proof of Lemma \ref{lem:suprho}, that
\beqn
\label{eq:2ndterm}
\E\abs{\rho_y-\hat\rho_y} \le \sqrt{\var[\hat\rho_y]}
\le
\sqrt{\rho_y\frac{1+2G\tha}{n(1-\tha)}}
,\qquad y\in\N.
\eeqn
To bound the first term, we recall the indicator variables $\xiy_i$ defined in (\ref{eq:xidef})
and observe that
\beq
n\E\abs{\rho_y-\hat\rho_y} &=& \E\abs{\sum_{i=1}^n\xiy_i-n\rho_y} \\
&\le& n\E\abs{\xiy_i-\rho_y} \\
&=& 2n\rho_y(1-\rho_y)\le 2n\rho_y,
\eeq
where stationarity was used in the last line of the derivation.
Combining the last display with (\ref{eq:2terms}) and (\ref{eq:2ndterm}) yields the claim.
\enpf
\bepf[Proof of Theorem \ref{thm:u-chernoff}.]
\bit
\item[(a)]
Since obviously
\beq
\narsum{\rho_y<1/n}\rho_y &\ninf& 0,
\eeq
it suffices to show that
\beqn
\label{eq:sqrtnp}
\oo{\sqrt{n}}\narsum{\rho_y\ge1/n}\sqrt{\rho_y} &\ninf& 0.
\eeqn
The latter was proved in \cite[Lemma 7]{Berend2013},
but we will present a simpler proof\footnote{
This elegant proof is due to Asaf Shachar. Andrew Barron points out that
(\ref{eq:sqrtnp}) may be easily derived from Lebesgue's dominated convergence theorem.}
 here.
Assume without loss of generality that $\rho_1\ge\rho_2\ge\ldots$, pick an abitrary $\eps>0$,
and let $N\in\N$ be large enough so that $\sum_{j\ge N}\rho_j<\eps$.
Then
\beq
\oo{\sqrt{n}}
\narsum{\rho_j\ge1/n}\sqrt{\rho_j}
&\le&
\oo{\sqrt{n}}\sum_{j\le N}\sqrt{\rho_j}
+
\oo{\sqrt{n}}\sum_{j>N,\rho_j\ge1/n}\sqrt{\rho_j} \\
&\le&
\sqrt{\frac{N}{n}}
+
\oo{\sqrt{n}}
\sqrt{\sum_{
\rho_j\ge1/n}1}
\sqrt{\sum_{j>N
}
\rho_j
}
\\
&\le&
\sqrt{\frac{N}{n}}+\sqrt{\eps},
\eeq
since there can be at most $n$ terms with $\rho_j\ge1/n$.
\item[(b)]
The claim follows from (\ref{eq:conc-tv}) and the fact that for any two distributions
$\phi,\psi\in\R^\N$,
\beq
\label{eq:TVE}
\TV{\phi-\psi}=\sup_{E\subseteq\N}\abs{\phi(E)-\psi(E)}.
\eeq
\eit
\enpf
\subsection{The stationarity assumption}
\label{sec:stat}
For rapidly mixing Markov and hidden Markov chains, the stationarity assumption can easily be relaxed.
Indeed, Let $Y=(Y_1,\ldots,Y_n)$ be a $(G,\tha)$-geometrically ergodic hidden Markov chain
with parameters $(B\pi',A,B)$,
where
$\pi'\in\R^\N$ is some
stochastic vector.
If $Y$
is ``nearly stationary,'' in the sense that
$\TV{\pi-\pi'}$ is small, 
a simple
dimension-free
bound 
on
the statistical distance between $Y$ and its stationary version is available.
\bethn
\label{thm:Y-Y'-dim-free}
Let $Y'=(Y'_1,\ldots,Y'_n)$ be the stationary version of $Y$ --- i.e., an HMM with parameters $(B\pi,A,B)$,
where $\pi$ is the stationary distribution of the kernel $A$.
Then
\beq
\TV{\calL(Y)-\calL(Y')}
\le \TV{\pi-\pi'}.
\eeq
\enthn
First, we prove an analogous result for Markov chains.
\belen
\label{lem:xi-xi'}
Let $A$ be 
Markov kernel 
and
$\xi,\xi'\in\R^\N$ two arbitrary
stochastic vectors. Let $X=(X_1,\ldots,X_n)$
and
$X'=(X'_1,\ldots,X'_n)$ be the Markov chains induced by $(\xi,A)$ and
$(\xi',A)$, respectively.
Then
\beq
\TV{\calL(X)-\calL(X')}=\TV{\xi-\xi'}
.
\eeq
\enlen
\bepf
\beq
\TV{\calL(X)-\calL(X')}
&=&
\oo2
\sum_{x\in\N^n}
\abs{
\paren{\xi_{x_1}-\xi'_{x_1}}
A_{x_2,x_1}\ldots A_{x_n,x_{n-1}}
}\\
&=&
\oo2
\sum_{x\in\N^n}
A_{x_2,x_1}\ldots A_{x_n,x_{n-1}}
\abs{
{\xi_{x_1}-\xi'_{x_1}}
}\\
&=&
\oo2
\sum_{x_1\in\N}
\abs{
{\xi_{x_1}-\xi'_{x_1}}
}=\TV{\xi-\xi'}.
\eeq
\enpf
\bepf[Proof of Theorem \ref{thm:Y-Y'-dim-free}]
Lemma \ref{lem:hmm<markov} lets us restrict our attention to the underlying Markov chains
$X$ and $X'$, respectively:
\beq
\TV{\calL(Y_{1\le i\le n})-\calL(Y'_{1\le i\le n})}
&\le&
\TV{\calL(X_{1\le i\le n})-\calL(X'_{1\le i\le n})}\\
&=&
\TV{\calL(X_{1})-\calL(X'_{1})}
= \TV{\pi-\pi'},
\eeq
where the first identity follows from Lemma \ref{lem:xi-xi'}.
\enpf
\becon
Let $Y_1,Y_2,\ldots$ be 
a 
(not necessarily stationary)
$\N$-valued
$(G,\theta)$-geometrically ergodic 
hidden Markov 
chain 
with 
stationary distribution 
$\rho=B\pi$ and initial distribution $\rho'=B\pi$.
Then the deviation bounds stated in Theorems \ref{thm:dkw} and \ref{thm:u-chernoff} hold
with an additive correction of $\TV{\pi-\pi'}$ on the right-hand side.
\encon
Letting a $(G,\tha)$-geometrically ergodic 
chain run for $s$ steps before starting the estimation ensures that $\TV{\pi-\pi'}\le G\tha^s$.
\subsection{Auxiliary lemma}
\label{sec:aux}
The Hamming metric 
on $\N^n$ 
is defined by 
$d(\bx,\by) = \sum_{i=1}^n \pred{x_i\neq y_i}$
for
$\bx,\by\in\N^n$.
\belen
\label{lem:Flip}
Suppose $n\in\N$ and $\bp\in\R^\N$ is 
a
distribution.
Define the functions $g,h:\N^n\to\R$ by
\beq
g(\bx) &=& \sup_{j\in\N}\abs{np_j-\sum_{i=1}^n\pred{x_i=j}},\qquad \bx\in\N^n,\\
h(\bx) &=& \sum_{j\in\N}\abs{np_j-\sum_{i=1}^n\pred{x_i=j}},\qquad \bx\in\N^n.
\eeq
Then 
$\Lip{g}\le1$ and
$\Lip{h}\le2$ with respect to the Hamming metric:
\beq
\abs{g(\bx)-g(\by)} &\le&  d(\bx,\by), \\
\abs{h(\bx)-h(\by)} &\le& 2 d(\bx,\by)
\eeq
for all $\bx,\by\in\N^n$.
\enlen
\bepf
We only prove the claim for $h$ (the proof for $g$ is analogous).
Let the function $\hat n_j:\N^n\to\N$ count the number of times $j$ appears in $x$;
formally, $\hat n_j(x) = \sum_{i=1}^n\pred{x_i=j}$. Now suppose $\bx,\by\in\N^n$
differ only in coordinate $k$, with $x_k=a$ and $y_k=b$.
Then
\beq
{h(\bx)-h(\by)} &=& 
{\sum_{j\in\N}\abs{np_j-\hat n_j(\bx)} - \sum_{j\in\N}\abs{np_j-\hat n_j(\by)}}\\
&=& 
{
\paren{
\abs{np_a-\hat n_a(\bx)}+\abs{np_b-\hat n_b(\bx)}
}
-
\paren{
\abs{np_a-\hat n_a(\by)}+\abs{np_b-\hat n_b(\by)}
}
}\\
&=&
{
\paren{
\abs{np_a-\hat n_a(\bx)}+\abs{np_b-\hat n_b(\bx)}
}
-
\paren{
\abs{np_a-(\hat n_a(\bx)-1)}+\abs{np_b-(\hat n_b(\bx)+1)}
}
}\\
&\le&
\dsabs{
\abs{np_a-\hat n_a(\bx)}
-
\abs{np_a-(\hat n_a(\bx)-1)}
}
+
\dsabs{
\abs{np_b-\hat n_b(\bx)}
-
\abs{np_b-(\hat n_b(\bx)+1)}
}\\
&\le&2.
\eeq
\enpf
\section*{Acknowledgments}
We thank the anonymous referee for carefully reading the manuscript and offering helpful suggestions.
\bibliographystyle{plain}
\bibliography{../../mybib}
\end{document}